\documentclass[11pt]{amsart}
\usepackage{graphpap}
\usepackage{amssymb, color}
\usepackage{amsmath}
\usepackage{amsfonts}
\newtheorem{theoreme}{Theorem}[section]
\newtheorem{lemma}{Lemma}

\newtheorem{proposition}{Proposition}

\newtheorem{remarque}[proposition]{Remark}

\numberwithin{equation}{section}

\def\R{\mathbb R}
\def\C{\mathbb C}
\def\Z{\mathbb Z}
\def\N{\mathbb N}

\def\T{\mathbb T}

\begin{document}
\title[Blow up in several points for NLS on a bounded domain]{Blow up in several points for the nonlinear Schr\"odinger equation on a bounded domain }
\author[N. Godet]{Nicolas Godet}
\address{University of Cergy-Pontoise, Department of Mathematics, CNRS, UMR 8088, F-95000 Cergy-Pontoise}
\email{nicolas.godet@u-cergy.fr}

\begin{abstract}
 
Given $p$ points in a bounded domain of $\R^d$, with $d=2,3$, we show the existence of solutions of the $L^2$-critical focusing nonlinear Schr\"odinger equation blowing up exactly in these points. 
\end{abstract}

\maketitle

\section{Introduction}

We consider the $L^2$-critical focusing nonlinear Schr\"odinger equation posed on a bounded and regular domain $\Omega$ of $\R^d$ (with $d=2,3$):
\begin{equation} \label{schro}
\begin{array}{l}
  i \partial_t u+ \Delta u=-|u|^{4/d}u,  \quad (t,x) \in [0,T) \times \Omega. 
\end{array}
\end{equation}
We add an initial data and the Dirichlet boundary condition :
\begin{equation} \label{schro2}
  \left \{
\begin{array}{l}
 u(t,x)=0, \quad (t,x) \in [0,T) \times \partial \Omega, \\
 u(0,x)=u_0(x), \quad x \in \Omega.
\end{array}
\right .
\end{equation}
Thanks to the Sobolev embedding $H^2(\Omega)\hookrightarrow L^{\infty}(\Omega)$, one can show that the equation (\ref{schro}) is locally well-posed in the space $H^2(\Omega) \cap H^1_0(\Omega)$ : for every initial data \sloppy $u_0 \in H^2(\Omega) \cap H^1_0(\Omega)$, \fussy there exists a time $T \in (0, + \infty] $ and a unique function \sloppy $u \in C([0,T),H^2(\Omega) \cap H^1_0(\Omega))$ \fussy solution of (\ref{schro}) with initial data $u_0$. If $u$ is a solution of (\ref{schro}), the energy and the mass are conserved : for every $t \in [0,T)$
\begin{eqnarray*}
 E(t)&:=&\frac{1}{2} \| \nabla u(t) \|_{L^2}^2- \frac{d}{4+2d} \| u(t)\|_{L^{\frac{4+2d}{d}}}^{\frac{4+2d}{d}}=E(0) ,\\
 M(t)&:=& \|u(t)\|_{L^2}=M(0).
\end{eqnarray*}
Moreover, we have the following blow up criteria
\[ \textrm{if} \quad T < + \infty \quad \textrm{then} \quad \|u(t)\|_{H^2 } \to + \infty \quad \textrm{when} \quad t \to T .\]
Note that if $d=2$, the equation (\ref{schro}) is also well posed in $H^1_0(\Omega)$ and even globally well posed if the $L^2$-norm of the initial data is smaller than $\|Q\|_{L^2(\R^2)}$ (see \cite{anton,brez,vlad}).

\vspace{0.2cm}

Blow up solutions for the equation (\ref{schro}) has been extremely studied in $\R^n$ and we expect that some results remain valid in bounded domains or more generally in the setting of flat geometries. Among papers concerning the study of the nonlinear Schr\"odinger equation on a domain, we can mention \cite{bani,fibi,kavian,planchon}. 

\vspace{0.2cm}

 In \cite{merle1}, Merle shows that if $\Omega$ is the whole space $\R^d$ (without restriction on $d$) then given $p$ points in $\R^d$, there exists a solution of the focusing nonlinear Schr\"odinger equation with $L^2$-critical nonlinearity that blows up in the $p$ points. The aim of this paper is to show that this result is still true if $\R^d$ is replaced by a bounded and regular domain of $\R^d$ with $d=2,3$. We can not adapt the construction of  Merle because the proof crucially uses the dispersion estimate
\[
  \exists C >0, \forall t >0, \forall v \in L^1(\R^d), \|e^{it \Delta}v\|_{L^{\infty}(\R^d)} \leq \frac{C}{t^{d/2} }\|v\|_{L^1(\R^d)}
\]
 which turns out to be false if we replace $\R^d$ by a bounded domain. To prove our result, we use the perturbation method introduced in \cite{ogawa} and used in \cite{niko} to treat the case of a point in dimension $2$. Because of lack of regularity of the nonlinearity, the case of the dimension $3$ requires a change in the choice of the weighted space where we show the property of contraction. 

\vspace{0.2cm}

In the following, we denote by $Q$ the unique (\cite{kwong,wein}) radially symmetric and strictly positive solution of 
\[ - \Delta Q + Q = |Q|^{4/d} Q \]
and satisfying the exponential decay at infinity :
\[ \forall \alpha \in \N^n, \exists C_{\alpha}>0, \exists D_{\alpha}>0, \forall x \in \R^d, | \partial^{\alpha} Q(x) | \leq C_{\alpha} e^{- D_{\alpha} |x|}.
\]

\section{Statement and proof of the result}
Now we state the theorem of the paper.
\begin{theoreme} \label{thm}
 Let $\Omega$ be a bounded and regular domain of $\R^d$ (with $d=2,3$) and $x_1, \dots,x_p$ $p$ distinct points in $\Omega$. Let $\varphi_1, \dots, \varphi_p \in C^{\infty}_0(\Omega)$ with disjoint supports and such that $0 \leq \varphi_k \leq 1$ and $\varphi_k=1$ near $x_k$. Then 

\begin{enumerate}
 \item  There exists $\lambda_0 >0$ such that for every $\lambda \geq \lambda_0$, there exists a time $T_{\lambda}>0$ and \sloppy $u_{\lambda} \in C([0,T_{\lambda}),H^2(\Omega) \cap H^1_0(\Omega))$ \fussy such that the function $h_\lambda$ defined by
\[ h_{\lambda}(t,x)=\frac{1}{\lambda ^{d/2} (T_{\lambda}-t)^{d/2}} \sum_{k=1}^p e^{\frac{i(4-\lambda^2 |x-x_k|^2)}{4 \lambda ^2(T_\lambda-t)}} \varphi_k(x) Q \left(\frac{x-x_k}{\lambda (T_{\lambda} -t)}\right) + u_{\lambda}(t,x) \]
is a solution of (\ref{schro}). Moreover, 
\[ \exists \ \gamma >0, C >0, \forall \lambda \geq \lambda_0, \forall t \in [0, T_{\lambda}), \ \|u_{\lambda}(t) \|_{H^2(\Omega)} \leq C  e^{- \frac{\gamma}{\lambda(T_{\lambda}-t)}}. \]
\item  For $\lambda \geq \lambda_0$, the solution $h_{\lambda}$ blows up in $H^1$ at time $t=T_{\lambda}$ in the points $x_1,\dots,x_p$ with speed $(T_{\lambda}-t)^{-1} $. More precisely $h_{\lambda}$ verifies
\begin{align*} 
(\romannumeral 1) & \ \ \textrm{for $R>0$ small enough, for every\ } k, \quad \| h_{\lambda}(t) \|_{L^2( \overline{B} (x_k,R))} \xrightarrow[t \to T_{\lambda}]{} \|Q\|_{L^2(\mathbb R^d)}, \\
  (\romannumeral 2) & \ \ \textrm{for all\ }  t \in[0,T_{\lambda}), \quad \|h_{\lambda}(t)\|_{L^2(\Omega)} =\sqrt{p} \|Q\|_{L^2(\mathbb R^d)},  \\
  (\romannumeral 3) & \ \ |h_{\lambda}(t)|^2  \underset{t \to T_{\lambda}}{\longrightarrow}\|Q\|_{L^2(\mathbb R^d)}^2 \sum_{k=1}^p  \delta_{x_k} \quad \textrm{in the sense of measures}, \\
  (\romannumeral 4) &  \ \ \| \nabla h_{\lambda}(t)\|_{L^2(\Omega)} \underset{t \to T_{\lambda}}{\sim}    \frac{\sqrt {p} \|\nabla Q \|_{L^2(\mathbb R^d)} }{ \lambda (T_{\lambda}-t)  }.
\end{align*}
\end{enumerate}
\end{theoreme}
{\textbf{Scheme of the proof.}} If $\Omega= \R^d$, we know an explicit solution $u_i$ of (\ref{schro}) which blows up in $x_i$. Next we consider the function $u=\varphi_1 u_1 + \dots\varphi_p u_p$ where $\varphi_k$ is a cut-off function near $x_k$. Therefore, $u$ is a function which vanishes on the boundary and has the same behavior than $u_k$ near $x_k$ because of the cut-off functions. Thus, $u$ blows up in the points $x_1,\dots,x_p$. However, $u$ is not a solution of $(\ref{schro})$ but we shall show the existence of a rest $r$ such that $u+r$ is a solution of $(\ref{schro})$ and keeps the behavior of $u$ when $t$ tends to the blow up time. For this, we will impose that $r(t)$ tends to $0$ at the blow up time. To prove the existence of the rest $r$, we perform a fixed point argument in a suitable weighted space.
\begin{proof} For $T>0$ and $\lambda>0$, we introduce
\[ \left\{ \begin{array}{rcl}
 & & \displaystyle{r^k_{\lambda}(t,x)=\frac{1}{\lambda ^{d/2}(T-t)^{d/2}} e^{\frac{i(4-\lambda^2 |x-x_k|^2)}{4 \lambda ^2(T-t)}}  Q \left( \frac{x-x_k}{\lambda (T -t)}  \right)} ,\\
 & & \displaystyle{r_{\lambda}(t,x)=\sum_{k=1}^p \varphi_k(x) r^k_{\lambda}(t,x) .}
 \end{array} \right.
 \]

For every $k,r_{\lambda}^k$ is a solution of (\ref{schro}) on $\R^d$. This solution blows up in $H^1(\R^d)$ at time $T$  and in $x_k$. We seek a condition on $u_{\lambda}$ for which the function $h_{\lambda}:=r_{\lambda} +u_{\lambda}$ is a solution of (\ref{schro}) on $\Omega$. We have the equality
\[ i\partial_t h_{\lambda} + \Delta h_{\lambda} =\sum_{k=1}^p  \left( r_{\lambda} ^k \ \Delta \varphi_k +2 \nabla \varphi_k \cdot \nabla r_{\lambda}^k -\varphi_k |r_{\lambda}^k|^{4/d}  r_{\lambda}^k \right)+(i\partial_t + \Delta) u_{\lambda} . \] 
Thus, $h_{\lambda}$ is a solution of (\ref{schro}) if and only if $u_\lambda$ satisfies
\begin{align*}
 (i\partial _t + \Delta)u_\lambda&=-|r_\lambda  + u_{\lambda} |^{4/d} (r_\lambda  + u_{\lambda}) -\sum_{k=1}^p  \left(r_{\lambda} ^k \ \Delta \varphi_k +2 \nabla \varphi_k \cdot \nabla r_{\lambda}^k -\varphi_k |u_{\lambda}^k|^{4/d}  u_{\lambda}^k \right), \\
  &= S_0+S(u_\lambda)
\end{align*}
where we denote
\[ \left \{ \begin{array}{rcl}
 S_0&= & \displaystyle{- \left| \sum_{k=1}^p \varphi_k r_\lambda^k \right|^{4/d} \left( \sum_{k=1}^p \varphi_k r_\lambda^k \right) -\sum_{k=1}^p  \left(r_{\lambda} ^k \ \Delta \varphi_k +2 \nabla \varphi_k \cdot \nabla r_{\lambda}^k -\varphi_k |r_{\lambda}^k|^{4/d}  r_{\lambda}^k \right)    },\\
 S(u)& =& \displaystyle{\left| \sum_{k=1}^p \varphi_k r_\lambda^k \right|^{4/d}  \left(\sum_{k=1}^p \varphi_k r_\lambda^k \right) -\left| u + \sum_{k=1}^p \varphi_k r_\lambda^k   \right|^{4/d} \left(u + \sum_{k=1}^p \varphi_k r_\lambda^k \right).  }\\
\end{array} \right. 
\]
To not perturb the behavior of $r_{\lambda}$ at the blow up time, we impose the condition $u(t)$ tends to $0$ when $t$ tends to $T$ ; this leads to consider the integral formulation
\begin{equation} \label{fixe}
 u(t)=i \int_t^T e^{i(t-s) \Delta} \left( S_0(s)+S(u)(s) \right) ds.
\end{equation}
We introduce
\[ \left\{
\begin{array}{l}
  I_0(t)= i \displaystyle{\int_t^T e^{i(t-s) \Delta} S_0(s) ds} ,\\
  I(u)(t)=i \displaystyle{\int_t^T e^{i(t-s) \Delta} S(u (s)) ds}. 
\end{array}
\right.
\]
To estimate the terms $I_0$ and $I$ we begin with the following lemma.

\begin{lemma} \label{lemestimate}
There exists a constant $C>0$ such that for every $u,v \in H^{2}(\Omega)$, 
\begin{align*} 
 (\romannumeral 1)& \ \| |u|^{4/d} u-|v|^{4/d}v \|_{L^2(\Omega)} \leq C \|u-v\|_{L^2(\Omega)} \left( \|u\|_{L^{\infty} (\Omega)}+ \|v\|_{L^{\infty} (\Omega)} \right)^{4/d}, \\
  (\romannumeral 2)& \ \| |u|^{4/d} u - |v|^{4/d}v \|_{H^1(\Omega)} \leq C \|u-v\|_{H^2(\Omega)} \left( \|u\|_{H^2(\Omega)}+ \|u\|_{H^2(\Omega)} \right )^{4/d}, \\
  (\romannumeral 3)& \    \| |u|^{4/d} u \|_{H^2(\Omega)} \leq C \|u\|_{H^2(\Omega)}^{1+4/d} .
\end{align*}
\end{lemma}

\begin{proof}[Proof of Lemma \ref{lemestimate}]
(\romannumeral 1). We use the Taylor formula
\begin{equation} \label{taylor}
 f(u)-f(v)=(u-v) \int_0^1 \partial_z f(tu+(1-t)v) dt + (\overline{u}-\overline{v}) \int_0^1 \partial_{\overline{z}} f(tu+(1-t)v) dt
\end{equation}
with the complex function $f(z)=|z|^{4/d}z$. The computation of the derivatives of $f$ shows that
\[
 |\partial_z f(z) | + |\partial_{\overline{z}} f(z) |  \leq C |z|^{4/d}.
\]
We deduce from (\ref{taylor}) that for $u,v \in \C$
\[
 \left | |u|^{4/d}u-|v|^{4/d}v \right |\leq C |u-v| \left ( |u|+ |v| \right)^{4/d}.
\]
Then, we apply this inequality to the functions $u$ and $v$, integrate and the conclusion follows by using H\"older inequality. This prove the first point.  

(\romannumeral 2). First, using (\romannumeral 1) and the Sobolev embedding $H^2(\Omega)\hookrightarrow L^{\infty}(\Omega)$  we get
\begin{equation} \label{norm}
 \| |u|^{4/d} u-|v|^{4/d} v \|_{L^2} \leq  C \|u-v\|_{H^2} \left ( \|u\|_{H^2(\Omega)} +\|v\|_{H^2(\Omega)}  \right )^{4/d}.
\end{equation} 
Next a direct computation shows that 
\begin{equation}
  \partial_i(|u|^{4/d}u)  = \left(\frac{2}{d}+1 \right) ( \partial_i u) |u|^{4/d} + \frac{2}{d} (\partial_i \overline{u}) |u|^{4/d-2}u^2 .  \label{deri1}
\end{equation}
According to the last identity, we can split $\partial_i(|u|^{4/d}u)-\partial_i(|v|^{4/d}v)$ into two terms. These two terms are treated in the same way. Let us for instance treat the second. We may write
\[ (\partial_i \overline{u}) |u|^{4/d-2}u^2 -(\partial_i \overline{v}) |v|^{4/d-2}v^2=A_1+A_2 ,\]
where
\[
 \left \{ 
\begin{array}{rcl}
 A_1&=& \partial_i(\overline{u}-\overline{v}) \left ( |u|^{4/d-2}u^2 \right), \\
A_2&=& \partial_i \overline{v} \left( |u|^{4/d-2} u^2-  |v|^{4/d-2} v^2  \right).
\end{array}
\right .
\]
But by H\"older inequality and again the Sobolev embedding of $H^2$ into $L^{\infty}$,
\begin{align*}
 \|A_1\|_{L^2} &\leq  \|u-v\|_{H^1} \|u\|_{L^{\infty}} ^{4/d} \\
               & \leq  \|u-v\|_{H^2} \|u\|_{H^2} ^{4/d}.
\end{align*}
For $A_2$, we first use (\ref{taylor}) to get
\[
  \left | |u|^{4/d-2}u^2- |v|^{4/d-2}v^2 \right | \leq C |u-v| \left( |u|+|v|\right)^{4/d-1}.
\]
Hence
\begin{align*}
 \|A_2\|_{L^2} &\leq \|v\|_{H^1} \|u-v\|_{L^{\infty}} \left ( \|u\|_{L^{\infty}} +\|v\|_{L^{\infty}} \right)^{4/d-1} \\
               & \leq \|u-v\|_{H^2} \left ( \|u\|_{H^2}+ \|v\|_{H^2} \right )^{4/d}.
\end{align*}
Summing the estimates on $A_1, A_2$ and (\ref{norm}), we obtain the second point. 

(\romannumeral 3). The norm $\| \cdot \|_{H^2}$ is equivalent to the norm $\| \cdot \|_{L^2}+ \| \nabla^2 \cdot \|_{L^2}$. Using again a Sobolev embedding, we deduce that 
\begin{equation}  \label{eq1}
\||u|^{4/d}u\|_{L^2} \leq  \|u\|_{L^{\infty}}^{4/d+1} \leq C \|u\|_{H^2}^{4/d+1}.
\end{equation}
Moreover, deriving the relation (\ref{deri1}), we get
\begin{align}
 \partial_{ij} (|u|^{4/d}u)=& \left( \frac{2}{d}+1 \right) \left( \partial_{ij}u |u|^{4/d} + \frac{2}{d} \partial_iu \partial_ju |u|^{4/d-2} \overline{u} + \frac{2}{d} \partial_i \overline{u} \partial_ju|u|^{4/d-2}u \right) \nonumber \\
 &+ \frac{2}{d} \left( \left (\frac{2}{d}+1 \right) \partial_iu \partial_j \overline{u} |u|^{4/d-2}u+ \partial_{ij} \overline{u}  |u|^{4/d-2}u^2 + \left (\frac{2}{d} -1 \right ) \partial_i \overline{u} \partial_j \overline{u} |u|^{4/d-4}u^3 \right ). \label{derivative}
\end{align}
Then using an H\"older inequality on each term, we get
\[
 \| \partial_{ij}   ( |u|^{4/d}u ) \|_{L^2} \leq C \left( \| \nabla ^2 u \|_{L^2} \|u\|_{L^{\infty}}^{4/d} + \|\nabla u\|_{L^4}^2 \|u\|_{L^{\infty}}^{4/d-1} \right).
\]
By the embeddings $H^1 \hookrightarrow L^4$ and $H^2 \hookrightarrow L^{\infty}$, we may write
\begin{equation} \label{eq2}
 \| \partial_{ij}  ( |u|^{4/d}u ) \|_{L^2} \leq C \|u\|_{H^2}^{1+4/d}.
\end{equation}
Gathering (\ref{eq1}) and (\ref{eq2}), we obtain the third point.
\end{proof}

Since the supports of $\varphi_k$ are disjoints, $S_0(t)$ is zero near the points $x_k$. Therefore there exists $r>0$ such that
\[
 \forall t \in [0,T), \ \forall x \in \Omega \setminus \cup_{k=1}^p  \overline{B}(x_k,r),\quad  S_0(t,x)=0. 
\]
Using the exponential decay of the ground state $Q$ and its derivatives, we get the existence of $C,D$ such that for every $k, t, \lambda, \alpha$ with $|\alpha| \leq 3$ :
\[   \| \partial ^{\alpha} r_{\lambda}^k(t)  \| _{L^2(\mathbb R^2 \setminus  \overline{B}(x_k,r))} \leq C e^{- \frac{D}{\lambda (T-t)}}. \]
Therefore, Lemma \ref{lemestimate} (\romannumeral 3) and the structure of algebra of $H^2(\Omega\setminus \cup_{k=1}^p \overline{B}(x_k,r))$ allows us to write
\begin{align} 
\| S_0(t) \|_{H^2(\Omega)}&=\| S_0(t) \|_{H^2(\Omega \setminus \cup_{k=1}^p \overline{B}(x_k,r))} \leq C e^{-\frac{\delta}{\lambda (T-t)}} .\label{estimate1}
\end{align}
Now, we can introduce the following metric space
\begin{equation*}
    E_{T}=  \left \{ u \in L^{\infty} ([0,T),H^2 \cap H^1_0 ), \sup_{0 \leq t < T} \left( e^{\frac{\delta}{\lambda (T-t)}} \|u(t)\|_{L^2} \right ) + \sup_{0 \leq t < T} \left ( e^{\frac{\alpha \delta}{\lambda (T-t)}} \|  u(t)\|_{H^2}  \right) \leq 1 \right \},
\end{equation*}
equipped with the distance
\[ d(u,v) = \sup_{0 \leq t < T} \left(  e^{\frac{\delta}{\lambda (T-t)}} \|u(t)-v(t)\|_{L^2} \right ),
\]
where $\alpha$ is a real number such that $0< \alpha <  \min(1,\frac{4}{d}-1)$. We are going to perform the Banach fixed point argument in $E_T$ to prove that the map defined by 
\begin{equation}
 \Phi(u)(t) = i\int_t^T  e^{i(t-s) \Delta } \left( S_0(s)+S(u)(s) \right)ds
\end{equation}
has a fixed point. For this, we show that for $T$ small enough and $\lambda$ big enough, $\Phi$ sends $E_T$ into itself and is a contraction.

First, we prove that $(E_T,d)$ is a complete metric space. It suffices to show that $E_T$ is closed in the complete metric space $E=\{ u \in L^{\infty}([0,T),L^2), \exp(\frac{\delta}{\lambda (T-t)}) \|u(t)\|_{L^2} \leq 1 \}$ equipped with the distance $d$. Let $(u_n)$ be a sequence in $E_T$ tending to $u \in E$ for $d$. Since $v_n(t):=\exp( \frac{\alpha \delta}{ \lambda (T-t)}) u_n(t)$ is bounded in $L^{\infty}([0,T), H^2 \cap H^1_0)$, we can extract a subsequence tending to $v \in L^{\infty}([0,T),H^2 \cap H^1_0)$ for the weak*-topology. Then necessarily, by uniqueness of the limit in $\mathcal D'((0,T),H^{-2})$, $v(t)=\exp( \frac{\alpha \delta}{ \lambda (T-t)}) u(t)$. Using the lower semi-continuity of the norm in $L^{\infty}([0,T),H^2 \cap H^1_0)$, we have
\[ 
  \sup_{0\leq t<T} \left (e^{ \frac{\alpha \delta}{ \lambda (T-t)}} \|  u(t)\|_{H^2} \right) \leq \liminf_{n \to \infty} \sup_{0 \leq t<T} \left (  e^{ \frac{\alpha \delta}{ \lambda (T-t)}} \| u_n(t)\|_{H^2} \right ), 
\]
hence taking the limit inferior in the inequality
\[
 \sup_{0 \leq t < T} \left( e^{\frac{\delta}{\lambda (T-t)}} \|u_n(t)\|_{L^2} \right ) + \sup_{0 \leq t < T} \left ( e^{\frac{\alpha \delta}{\lambda (T-t)}} \|  u_n(t)\|_{H^2}  \right) \leq 1  ,
\]
we obtain that $u \in E_T$.

\vspace{0.2cm}

\indent {\textbf{Boundedness.}} Let $T \in (0,1)$, $\lambda \geq 1$ to be chosen later. We prove that $\Phi$ sends $E_T$ in itself.  We may write
\[
 S(u)-S(v) =|r_{\lambda}+v|^{4/d}(r_{\lambda}+v)-|r_{\lambda}+u|^{4/d}(r_{\lambda}+u).
\]
Lemma \ref{lemestimate} and the Sobolev embedding $H^2 \hookrightarrow L^{\infty}$ provides
\begin{align*}
  \|S(u)(t)-S(v)(t)\|_{L^2} & \leq C \|u(t)-v(t)\|_{L^2} \left( \|u(t)\|_{H^2} + \|v(t) \|_{H^2} + \|r_{\lambda}(t) \|_{L^\infty} \right)^{4/d} .
\end{align*}  
But using the explicit formula of $r_{\lambda}$, one can compute its derivatives to get for $k=0,1,2$,
\begin{equation}
 \|\nabla^k r_{ \lambda}(t) \|_{L^{\infty}(\Omega)} \leq \frac{C}{ \lambda ^{d/2} (T-t)^{d/2+k}}, \label{deri}
\end{equation}
and deduce that if $u, v \in E_T$,
\[ \|S(u)(t)-S(v)(t)\|_{L^2} \leq C e^{- \frac{\delta}{ \lambda (T-t)} } d(u,v) \left (       1 + \frac{1}{\lambda ^{2}(T-t)^{2}} \right) .
\]
By integrating this inequality, and using
\[ 
 \int_t^T \frac{e^{ - \frac{\delta}{\lambda (T-s)} } }{(T-s)^2} ds \leq C \lambda e^{ - \frac{\delta}{\lambda (T-t)}} ,
\]
it follows 
\begin{align} \label{ineg1}
 \|I(u)(t)-I(v)(t) \|_{L^2} & \leq C  e^{- \frac{\delta}{ \lambda (T-t)} } d(u,v) \left ( T  + \frac{1}{\lambda}  \right ) .
\end{align}

Now we need to bound $\|  I(u)(t)\|_{H^2}$. The norm $\| \cdot \|_{H^2}$ is equivalent to $\| \cdot \|_{L^2}+\| \nabla ^2 \cdot \|_{H^2} $. In the formula (\ref{derivative}), there are two types of terms. The first one are those containing the product of two first derivatives, namely
\[ \partial_iu \partial_ju |u|^{4/d-2} \overline{u}, \ \ \partial_i \overline{u} \partial_ju|u|^{4/d-2}u, \ \ \partial_iu \partial_j \overline{u} |u|^{4/d-2}u, \ \ \partial_i \overline{u} \partial_j \overline{u} |u|^{4/d-3} u^3 ; \]
and the second are those containing a second derivative, namely
\[ \partial_{ij} u |u|^{4/d}, \ \ \ \ \ \partial_{ij} \overline{u} |u|^{4/d-2} u^2.\]
Each term of the same type are treated in the same way. Let us begin with the first type, for instance the term $\partial_i\overline{u} \partial_j\overline{u} |u|^{4/d-4} u^3$. We may write
\begin{align*}
\left \| \partial_i \overline{r_{\lambda}} \partial_j \overline{r_{\lambda}} |r_{\lambda}|^{4/d-4}r_{\lambda}^3 - \partial_i (\overline{r_{\lambda}}+\overline{u}) \partial_j (\overline{r_{\lambda}}+ \overline{u})  |r_{\lambda}+u|^{4/d-4}(r_{\lambda}+u)^3 \right \|_{L^2}  & \leq B_1+B_2
\end{align*}
where
\[
 \left \{ 
\begin{array}{rcl}
 B_1&=& \displaystyle{ \left \| \partial_i \overline{r_{\lambda}} \partial_j \overline{r_{\lambda}} \left( |r_{\lambda}|^{4/d-4} r_{\lambda}^3-|r_{\lambda}+u|^{4/d-4}(r_{\lambda} + u)^3 \right) \right \|_{L^2} ,}\\
 B_2&=& \displaystyle{\left \|  \left(\partial_i \overline{r_{\lambda}} \partial_j \overline{u}+ \partial_i \overline{u} \partial_j \overline{r_{\lambda}}+ \partial_i \overline{u} \partial_j \overline{u}\right) |r_{\lambda}+u|^{4/d-4}(r_{\lambda} + u)^3 \right \|_{L^2}.}
\end{array}
\right.
\]
Using the estimate on the derivatives of $r_{\lambda}$ (\ref{deri}) and the inequality 
\[ \left |  |u|^{4/d-4}u^3- |v|^{4/d-4}v^3 \right | \leq C |u-v|^{4/d-1},
\]
we have successively (with always $u \in E_T$).
\begin{align*}
 B_1 & \leq C  \| \partial_i r_{\lambda} \partial_j r_{\lambda} |u|^{4/d-1} \|_{L^2} \\
   & \leq C \| \nabla r_{\lambda} \|^2_{L^{\infty}} \|  |u|^{4/d-1}\|_{L^2} \\
   & \leq  \frac{C}{ \lambda^{d} (T-t)^{d+2}} \|u\|_{L^2}^{4/d-1} \\
   & \leq C \lambda^2 e^{-\frac{\alpha \delta}{ \lambda (T-t)} }. 
\end{align*}
For $B_2$, H\"older inequality provides
\begin{align*}
 B_2 &\leq  \| |r_{\lambda}+u|^{4/d-1} \left( \partial_i r_{\lambda} \partial_j u+ \partial_i u \partial_j r_{\lambda}+\partial_i u \partial_j u \right) \|_{L^2} \\
& \leq \|r_{\lambda}+u \|_{L^{\infty}}^{4/d-1} \left(  \| \nabla r_{\lambda} \|_{L^{\infty}} \| \nabla u \|_{L^2}+ \|\nabla u\|_{L^4}^2 \right) .
\end{align*}
But by interpolation between the spaces $L^2$ and $H^2$, we remark that there exists $\beta \in (\alpha,1)$ such that for every $u \in E_T$,  and for every $t \in [0,T)$,
\[
 \|u(t) \|_{H^1} \leq e^{- \frac{\beta \delta}{ \lambda (T-t)}} .
\]
Hence
\begin{align*}
B_2 & \leq  C \left( \frac{1}{\lambda ^{2-d/2}(T-t)^{2-d/2}} + 1 \right) \left( \frac{1}{\lambda^{d/2}(T-t)^{d/2+1}} e^{-\frac{ \beta  \delta}{ \lambda (T-t)}}  +  e^{-\frac{2 \alpha \delta}{ \lambda (T-t)}}   \right) \\
& \leq C \lambda e^{-\frac{\alpha \delta}{ \lambda (T-t)}} .
\end{align*}
Now we treat terms belonging to the second type. For instance the term $\partial_{ij} \overline{u} |u|^{4/d-2} u^2$. We have
\[ \| \partial_{ij} \overline{r_{\lambda}} |r_{\lambda}|^{4/d-2} r_{\lambda}^2-\partial_{ij} (\overline{r_{\lambda}}+ \overline{u}) |r_{\lambda}+u|^{4/d-2} (r_{\lambda}+u)^2 \|_{L^2} \leq C_1+C_2,\]
where
\[
 \left \{ 
\begin{array}{rcl}
 C_1&=& \displaystyle{ \left \| \partial_{ij} \overline{r_{\lambda}} \left( |r_{\lambda}|^{4/d-2} r_{\lambda}^2-|r_{\lambda}+u|^{4/d-2} (r_{\lambda}+u)^2 \right) \right \|_{L^2},} \\
 C_2&=& \displaystyle{\left \| \partial_{ij} \overline{u} \left ( |r_{\lambda}+u|^{4/d-2} (r_{\lambda}+u)^2\right) \right \|_{L^2} }.
\end{array}
\right.
\]
The inequality
\[
\left | |z|^{4/d-2}z^2-|w|^{4/d-2}w^2 \right | \leq C |z-w| \left (  |z|+ |u| \right)^{4/d-1},
\]
applied with $z=r_{\lambda}$ and $w=r_{\lambda}+u$ allows to write
\begin{align*}
C_1 & \leq \|\nabla^2 r_{\lambda} \|_{L^{\infty}} \| u \|_{L^2}  \left( \|r_{\lambda} \|_{L^ \infty}^{4/d-1} + \|u\|_{L^{ \infty}}^{4/d-1} \right) \\
    & \leq \frac{C}{ \lambda^{d/2} (T-t)^{d/2+2} } e^{- \frac{\delta}{ \lambda (T-t)}} \left( \frac{1}{ \lambda^{2-d/2} (T-t)^{2-d/2}}+1  \right)  \\
    & \leq C \lambda ^2 e^{- \frac{\alpha \delta}{  \lambda (T-t)} } .
\end{align*}
For the second bound, we may write
\begin{align*}
C_2 & \leq \| \nabla^2 u \|_{L^{2}}  \|r_{\lambda}+u\|_{L^{\infty}}^{4/d} \\
    & \leq C e^{- \frac{\alpha \delta}{   \lambda (T-t)} } \left (  \frac{1}{ \lambda^2 (T-t)^2}  +1\right).
\end{align*}
Summing $B_1, B_2, C_1$ and $C_2$, we deduce for $\lambda \geq 1$ and $u \in E_T$,
\begin{equation} \label{estiS}
 \|S(u)(t)\|_{H^2} \leq C e^{- \frac{\alpha \delta}{  \lambda (T-t)}} \left ( \lambda ^2 + \frac{1}{ \lambda^2 (T-t)^2} \right)
\end{equation}
and integrating we get 
\begin{align} 
\|  I(u)(t) \|_{H^2} & \leq C e^{- \frac{\alpha \delta}{   \lambda (T-t)} } \left( T \lambda ^2 + \frac{1}{\lambda} \right) .\label{estiI}
\end{align}
Moreover by (\ref{estimate1}), we obtain the bound
\[ \|I_0\|_{H^2} \leq C T e^{- \frac{\delta}{\lambda (T-t)}}.\]
Finally, the latter estimate together with (\ref{estiI}) and (\ref{ineg1}) applied with $v=0$ gives for $\lambda \geq 1$ and $u \in E_T$
\begin{equation} \label{estimate3}
 \sup_{0 \leq t <T} \left ( e^{ \frac{\delta}{\lambda (T-t)}} \|\Phi(u)(t)\|_{L^2} \right ) + \sup_{0 \leq t <T} \left ( e^{ \frac{\alpha \delta}{\lambda (T-t)}} \|\Phi (u)(t)\|_{H^2} \right ) \leq C \left ( T \lambda ^2 + \frac{1}{\lambda} \right ).
\end{equation}

\vspace{0.2cm}

\textbf{Contraction.} Actually, we have already proved the property of contraction of $\Phi$ during the proof of the boundedness. Indeed, the estimate (\ref{ineg1}) provides for every $u,v \in E_T$,
\begin{equation} \label{estifi}
 d \left( \Phi(u),\Phi(v) \right) \leq C \left (T + \frac{1}{\lambda} \right ) d(u,v).
\end{equation}

\textbf{Conclusion.} First $\Phi(E_T) \subset L^{\infty}([0,T), H^2 \cap H^1_0)$. Indeed, the bound (\ref{estimate3}) gives $\Phi( E_T) \subset L^{\infty}([0,T), H^2)$ and it remains to verify that $\Phi(u)(t) \in H^1_0$ almost everywhere. For this, it suffices to prove that for $u \in H^2 \cap H^1_0, S(u) \in H^1_0$. But if $u \in H^2 \cap H^1_0$, we can appoximate $u$ in $H^2$-norm by a sequence $u_n \in C^{\infty}_0$. And by Lemma \ref{lemestimate} (\romannumeral 2), we obtain the convergence in $H^1$-norm of $|u_n|^{4/d}u_n$ to $|u|^{4/d}u$ and this shows that $|u|^{4/d}u \in H_0^1$. Hence $S(u) \in H^1_0$ and $\Phi(E_T) \subset L^{\infty}([0,T), H^2 \cap H^1_0)$. Moreover estimates (\ref{estimate3}) and (\ref{estifi}) prove that we can choose $\lambda$ big enough and $T>0$, depending on $\lambda$,  such that $\Phi( E_T) \subset E_T$ and for every $u,v \in E_T$,
\begin{equation}  \label{final}
 d \left (  \Phi(u),\Phi(v) \right )  \leq  \frac{1}{2} d \left( u,v \right) .
\end{equation}

Thus, we can apply the Banach fixed point argument with the function $\Phi$ and this proves the existence of the rest $u \in E_T$ satisfying (\ref{fixe}). To obtain the continuity in time with values in $H^2$, we use the integral formulation satisfied by $u$. Estimates (\ref{estimate1}) and (\ref{estiS}) show that the map $s \mapsto S_0(s)+(S(u))(s)$ belongs to $L^{\infty}([0,T),H^2)$ and since $u$ verifies
\[ u(t)=i e^{it \Delta} \int_t^T e^{-i s \Delta} \left( S_0(s) + (S(u))(s) \right) ds ,\]
we conclude that $u \in C([0,T),H^2)$.

\vspace{0.2cm}

Let us show the second part of the theorem. (\romannumeral 1). For $R>0$ small enough, we remark that for $x \in \cup_k \overline{B}(x_k,R)$
\[ r_{\lambda}(t,x)=\sum_{k=1}^p r^k_{\lambda}(t,x) . \]
Then
\[ \left | \|h_{\lambda}(t)\|_{L^2(\overline {B}(x_k,R))}  - \|r_{\lambda}^k(t)\|_{L^2(\overline {B}(x_k,R))} \right| \leq \|u_{\lambda}(t)\|_{L^2(\overline {B}(x_k,R))} + \sum_{j \neq k} \|r_{\lambda}^j(t)\|_{L^2(\overline {B}(x_k,R))} . \]
Using the fact that the $L^2$ norm of $r_{\lambda}^k$ is only concentrated in $x_k$ with mass $\|Q\|_{L^2(\R^d)}$, and that $\|u_{\lambda}(t)\|_{L^2}$ converges to $0$ when $t$ goes to $T_{\lambda}$, we obtain $(\romannumeral 1)$ by a passage to the limit. 

(\romannumeral 2). We may write
\[
 \| h_{\lambda}(t)\|_{L^2}^2= \|u_{\lambda }(t)\|_{L^2}^2 +\sum_{k=1}^p \| \varphi_k r_{\lambda}^k(t) \|_{L^2(\Omega)}^2 +2 \mathrm{Re} ( \langle r_{\lambda}(t),u_{\lambda}(t) \rangle _{L^2} ).
\]
The first and the latter term goes to $0$ when $t$ tends to $T_{\lambda}$ because $\|u_{\lambda} (t) \|_{L^2}$ tends to $0$. And the second term converges to $p \|Q\|_{L^2}^2$ thanks to the property of concentration of $r_{\lambda}^i$ near $x_i$. 

(\romannumeral 3). Let $\psi $ be a continuous function with compact support. Then if we denote
\[ I(t)=\int_{\Omega} | h_{\lambda}(t,x)|^2 \psi(x) dx -\|Q\|_{L^2(\mathbb R^d)}^2 \sum_{k=1}^p \psi(x_k)\ , \]
we have
\begin{align*}
\left | I(t) \right| & \leq   C \int_{\mathbb R^d} \left| | r_{\lambda}(t,x)|^2 \psi(x)  -Q^2(x) \sum_{k=1}^p \psi(x_k) \right| dx + C \int_{\Omega} | u_{\lambda}(t,x) |^2 dx  \\
                                            &  \leq C \sum_{k=1}^p \int_{\mathbb R^d}  \left | \left|r_{\lambda}^k(t,x) \right|^2 \psi(x)-Q^2(x) \psi(x_k) \right| dx + Ce^{- \frac{2 \delta}{\lambda (T_{\lambda}-t)}}   \\
                                            &\mbox{} \quad +  C \sum_{k=1}^p \int_{\mathbb R^d \setminus \overline{B}(x_k,R)} \left| (1-\varphi_k^2(x)) \psi(x) \left|r_{\lambda}^k(t,x) \right|^2  \right| dx \\
                                            & \leq C \sum_{k=1}^p \int_{\mathbb R^d}  \left | \left|r_{\lambda}^k(t,x) \right|^2 \psi(x)-Q^2(x) \psi(x_k) \right| dx + C \sum_{k=1}^p \int_{\mathbb R^d \setminus \overline{B}(x_k,R)}  \left|r_{\lambda}^k(t,x) \right|^2   dx \\
                                            &  \mbox{} \quad + C e^{- \frac{2 \delta}{\lambda (T_{\lambda}-t)}}  .                             
\end{align*}
But $|r_{\lambda}^k(t)|^2$ converges to $\|Q\|_{L^2(\R^d)}^2 \delta_{x_k}$ when $t \to T_{\lambda}$ so the first term goes to $0$. And the second one as well by the well known properties of $r_{\lambda}^k$. 

(\romannumeral 4). We have the equality
\[ \nabla h_{\lambda}= \nabla u_{\lambda} + \sum_{k=1}^p  r_{\lambda}^k \nabla \varphi_k + \sum_{k=1}^p \varphi_k \nabla r_{\lambda}^k . \]
Remarking that $\|\nabla u_{\lambda}(t) \|_{L^2(\Omega)}$ decays to $0$ when $t$ goes to $T_{\lambda}$ and  $\|\sum_{k=1}^p r_{\lambda}^k(t) \nabla \varphi_k \|   _{L^2(\Omega)}$ is bounded, we get the equivalence
\begin{align*}
\| \nabla h_{\lambda}(t) \|_{L^2(\Omega)} & \underset{t \to T_{\lambda}}{\sim} \left \| \sum_{k=1}^p  \varphi_k \nabla r_{\lambda}^k(t)   \right \|_{L^2(\Omega)}.
\end{align*}
But the $\varphi_k$ have disjoint supports so
\begin{align*}
\left \| \sum_{k=1}^p  \varphi_k \nabla r_{\lambda}^k(t)   \right \|_{L^2(\Omega)}&=  \left( \sum_{k=1}^p \left \| \varphi_k \nabla r_{\lambda}^k(t)   \right \|_{L^2(\Omega)}^2 \right)^{1/2} ;
\end{align*} 
and for all $k$,
\begin{align*}
\left \| \varphi_k \nabla r_{\lambda}^k(t)   \right \|_{L^2(\Omega)} & \underset{t \to T_{\lambda}}{\sim} \left \| \nabla r_{\lambda}^k(t)   \right \|_{L^2(\Omega)} \\ 
                            & \underset{t \to T_{\lambda}}{\sim} \frac{1}{\lambda (T_{\lambda}-t)} \|\nabla Q \|_{L^2(\mathbb R^d)}.
\end{align*}
We obtain (\romannumeral 4) by summing these equivalences.
\end{proof}

Let us give some remarks about Theorem \ref{thm}.
\begin{remarque}
The existence of blow up solutions for the $L^2$-critical focusing nonlinear Schr\"odinger equation remains true if we replace the bounded domain of $\R^d$ by a Riemannian manifold of dimension $d=2,3$ which is locally isometric to an open subset of $\R^d$ near the blow up points ; that is the case  of the flat torus $\T^d=\R^d / \Z^d$.
\end{remarque}
\begin{remarque}
For $d=3,4$, one can construct solutions of the following $L^2$-supercritical equation posed on $\T^{d}$
\begin{equation} \label{eq}
i \partial_t u + \Delta_{\T^{d}} u = - |u|^{\frac{4}{d-1}}u , 
\end{equation}
blowing up on the union the $p$ circles. Indeed, for $x_1,\dots,x_p \in \T^{d-1}$, by Theorem \ref{thm}, there exists a solution $u \in C([0,T),H^2(\T^{d-1}))$ blowing up in the $p$ points. Then, we consider the function $v \in C([0,T), H^2(\T^{d}))$ defined by $v(t,a,b)=u(t,a)$ for $a \in \T^{d-1}, b \in \T$. Then $v$ is a solution of (\ref{eq}) and blows up in the $p$ circles $\{x_k \} \times \T$. Note that the blow up on a sphere for the supercritical equation has been studied more precisely in $\R^n$ (see \cite{raphael}).
\end{remarque}

\begin{remarque}
 One can also interest in the case where the equation is posed in dimension greater than $3$. However, in this case, the nonlinearity is not regular enough to perform the same proof than in the case $d \leq 3$. Indeed, we need to solve the equation in a space included in $H^s$ with $s> d/2 \geq 2$ to get the embedding into $L^{\infty}$ but if $d \geq 4$, we can not derive the nonlinearity more than two times.
\end{remarque}

\noindent \textbf{Acknowledgements.} I would like to thank Nikolay Tzvetkov for introducing this subject, for his advices and his helpful remarks on my work.


\begin{thebibliography}{13}
\bibitem{anton} R. ANTON, \textit{Strichartz inequalities for Lipschitz metrics on manifolds and nonlinear Schr\"odinger equation on domains}, Bull. Soc. Math. Fr.
\bibitem{bani} V. BANICA, \textit{Remarks on the blow-up for the Schr\"odinger equation with critical mass on a plane domain.}, Ann. Sc. Norm. Super. Pisa (5), Vol. III (2004), 139-170.
\bibitem{brez} H. BREZIS, T. GALLOUET, \textit{Nonlinear Schr\"odinger evolution equation, Nonlinear Analysis}, Theory Methods Appl. 4 (1980), 677-681.
\bibitem{niko}  N. BURQ, P. G\'{E}RARD, N. TZVETKOV, \textit{Two singular dynamics of the nonlinear Schr\"odinger equation on a plane domain}, Geom. Funct. Anal. 13 (2003), 1-19.
\bibitem{caze} T. CAZENAVE, \textit{An introduction to nonlinear Schr\"odinger equations}, Text. Met. Mat 22, Inst. Math, Rio de Janeiro, (1989).
\bibitem{fibi} G. FIBICH, F. MERLE, \textit{Self-focusing on bounded domains}, Phys. D 155 (2001), no. 1-2, 132--158.
\bibitem{kavian} O. KAVIAN, \textit{A remark on the blowing up of solutions to the Cauchy problem for nonlinear Schr\"odinger equations}, Trans. AMS 299 (1987), 193-203.
\bibitem{kwong} M.K. KWONG, \textit{Uniqueness of positives solutions of $\Delta u-u + u^p=0$ in $\R^N$}, Arch. Rat. Mech. Ann. 105 (1989), 243-266.
\bibitem{merle1} F. MERLE \textit{Construction of solutions with exactly $k$ Blow-up points for the Schr\"odinger equation with critical nonlinearity}, Commun. Math. Phys. 129, 223-240 (1990).
\bibitem{ogawa} T. OGAWA, Y. TSUTSUMI, \textit{Blow-up solutions for the nonlinear Schr\"odinger equation with quartic potential and periodic boundary conditions}, Springer Lectures Notes in Math. 1450 (1990), 236-251.
\bibitem{planchon} F. PLANCHON, P. RAPHAEL, \textit{Existence and stability of the log-log blow up dynamics for the $L^2$ critical nonlinear Schr\"odinger equation in a domain}, Ann. Henri Poincar\'e 8 (2007), no. 6, 1177-1219.
\bibitem{raphael} P. RAPHAEL, \textit{Existence and stability of a solution blowing up on a sphere for an $L^2$-supercritical nonlinear Schr\"odinger equation}, Duke Math. J. Volume 134, Number 2 (2006), 199-258.
\bibitem{vlad} M. V. VLADIMIROV, \textit{On the solvability of mixed problem for a nonlinear equation of Schr\"odinger type}, Dokl. Akad. Nauk SSSR 275 (1984), 780-783.
\bibitem{wein} M.I. WEINSTEIN, \textit{Nonlinear Schr\"odinger equations and sharp interpolation estimate}, Commun. Math. Phys. 87 (1983) 567-576.


 
                               

\end{thebibliography}
\end{document}